\documentclass[10pt, a4paper]{amsart}
\usepackage{amssymb,latexsym,amsxtra,amscd,amsfonts,amsthm,amsmath,euscript}

\def\Z{{\mathbb Z}}
\def\NN{{\mathbb N}}
\def\PP{{\mathbb P}}

\def\Fq{{\mathbb F}_q}
\def\LA{{\mathcal L}_A}

\def\a{{\alpha}}
\def\b{{\beta}}
\newcommand\codim{\mathop{\rm codim}}
\newtheorem{theorem}{Theorem}
\newtheorem{proposition}[theorem]{Proposition}
\newtheorem{lemma}[theorem]{Lemma}
\newtheorem{corollary}[theorem]{Corollary}

\newtheorem*{theorem*}{Theorem}
\newtheorem*{conjecture*}{Conjecture}
\theoremstyle{remark}
\newtheorem{remark}[theorem]{Remark}

\begin{document}
\title{Schubert Varieties, Linear Codes and Enumerative Combinatorics}

\author{Sudhir R. Ghorpade}
\thanks{The first named author was partially supported by a `Career Award' grant
from AICTE, New Delhi and an IRCC grant from IIT Bombay.
The second named author was partially supported by 
the RFBR Grants 99-01-01204, 02-01-01041 and 02-01-22005.}
\address{Department of Mathematics, 
Indian Institute of Technology Bombay,\newline \indent
Powai, Mumbai 400076, India.}
\email{srg@math.iitb.ac.in}

\author{Michael A. Tsfasman}

\address{
Institut de Math\'ematiques de Luminy,
Case 907, 13288 Marseille, France,\newline \indent
and \newline \indent
Independent University of Moscow,\newline \indent
and \newline \indent
Dorbushin Math. Lab., Institute for Information Transmission Problems, Moscow.}
\email{tsfasman@iml.univ-mrs.fr}

\date{\today}

\keywords{Grassmannian, linear codes, minimum distance, projective system, 
Schubert variety} 
\begin{abstract}
We consider linear error correcting codes associated to higher dimensional
projective varieties defined over a finite field. The problem of
determining the basic parameters of such codes often leads to some
interesting and difficult questions in combinatorics and algebraic geometry. 
This is illustrated by codes associated to Schubert varieties in
Grassmannians, called  Schubert codes, which have recently been studied. 
The basic parameters such as the length, dimension and minimum distance 
of these codes are known only in special cases.  An upper bound for 
the minimum distance is known and it is conjectured that this bound is achieved.
We give explicit formulae for the length and dimension of arbitrary Schubert codes 
and prove the minimum distance conjecture in the affirmative for codes associated 
to Schubert divisors.
\end{abstract}

\maketitle
 
\section{Introduction}
\label{sec:intro}

Let $\Fq$ denote the finite field with $q$ elements, and let $n,k$ be integers
with $1\le k\le n$. 
The $n$-dimensional vector space $\Fq^n$ has a  norm, 
called {\em Hamming norm}, which is defined by
$$
\Vert x \Vert = \vert \left\{i\in\{1, \dots , n\} : x_i\ne 0\right\} \vert
\quad \mbox{for $x\in \Fq^n$.}
$$
More generally, if $D$ is a subspace of $\Fq^n$, the {\em Hamming norm
of $D$} is defined by 
$$
\Vert D \Vert = \vert \left\{i\in\{1, \dots , n\} : \mbox{ there exists $x\in D$
with } x_i\ne 0\right\} \vert .
$$

A {\em linear $[n,k]_q$-code} 
is, by definition, a $k$-dimensional subspace of $\Fq^n$. 
The adjective {\em linear} will often be dropped since in this paper we only consider
linear codes.
The parameters $n$ and $k$ are referred to as the {\em length} and the 
{\em dimension} of the corresponding code.
If $C$ is an $[n,k]_q$-code, then the {\em minimum distance} $d = d(C)$ of $C$
is defined by
$$
d(C) = \min\left\{ \Vert x \Vert : x\in C, \ x\ne 0\right\}.
$$
More generally, given any positive integer $r$, the {\em $r$th higher
weight} $d_r = d_r(C)$ of $C$ is defined by
$$
d_r(C) = \min\left\{ \Vert D \Vert : D \mbox{ is a subspace of $C$ with } \dim D =r\right\}.
$$
Note that $d_1(C) = d(C)$.

An $[n,k]_q$-code is said to be {\em nondegenerate} if it is 
not contained in a coordinate hyperplane of $\Fq^n$.
Two $[n,k]_q$-codes 
are said to be {\em equivalent} if one
can be obtained from another by permuting coordinates and multiplying them by
nonzero elements of $\Fq$; in other words, if they are in the same orbit
for the natural action of the semidirect product of $(\Fq^*)^n$ and
$S_n$. It is clear that this gives a natural equivalence relation on the set of
$[n,k]_q$-codes.

An alternative way to describe codes is via the language of projective systems
introduced in \cite{TV}.
A {\em projective system} is a (multi)set $X$ of $n$ points 
in the projective space $\PP^{k-1}$ 
over $\Fq$. We call
$X$ {\em nondegenerate} if these $n$ points are not contained in a
hyperplane of $\PP^{k-1}$. Two projective systems in $\PP^{k-1}$ are said to be 
{\em equivalent} if there is a projective automorphism of the ambient space
$\PP^{k-1}$, which maps one to the other; in other words, if they are in the 
same orbit for the natural action of $PGL(k, \Fq)$. It is clear that this 
gives a natural equivalence relation on the set of projective systems of 
$n$ points in $\PP^{k-1}$.

It turns out that a nondegenerate projective system of $n$ points in 
$\PP^{k-1}$ corresponds naturally to a nondegenerate linear $[n,k]_q$-code. 
Moreover, if we pass to equivalence classes with respect to the equivalence 
relations defined above, then this correspondence is one-to-one.
The minimum distance of the code $C=C_X$ associated to a 
nondegenerate projective system $X$ of $n$ points in 
$\PP^{k-1}$ admits a nice geometric interpretation in terms of $X$, namely, 
$$
d(C_X)  = 
n - \max\left\{|X\cap H| : H \mbox{ a hyperplane of } \PP^{k-1} \right\}. \\
$$
We have a similar interpretation for the $r$th higher weight $d_r(C_X)$, where the 
hyperplane $H$ is replaced by a projective subspace of codimension $r$ in $\PP^{k-1}$.
For more details concerning projective systems, higher weights and a proof of the
above mentioned one-to-one correspondence, we refer to \cite{TV} and \cite{TV2}.

The language of projective systems
not only explains the close connection between algebraic geometry and coding
theory, but also facilitates the introduction of linear codes corresponding to
projective algebraic varieties defined over a finite field.
A case in point is the Grassmannian $G_{\ell, m} = G_{\ell}(V)$ of 
$\ell$-dimensional subspaces of an $m$-dimensional vector space $V$ over 
$\Fq$. We have the well-known
Pl\"ucker embedding of the Grassmannian into a projective space (cf. \cite{GL}, 
\cite{HP}), and this embedding is known to be nondegenerate. 
Considering the ($\Fq$-rational) points of  $G_{\ell, m}$ as a projective system,
we obtain a $q$-ary linear code, called the {\em Grassmann code}, which 
we denote by $C(\ell ,m)$. These
codes were first studied by Ryan \cite{CR1,CR2, CKR} in the binary case
and by  Nogin \cite{N} in the $q$-ary case.  It is 
clear that the length $n$ and the dimension $k$ of $C(\ell ,m)$ are given by
\begin{equation}
\label{eq:1}
n = { {m} \brack {\ell} }_q :=
 \frac{ (q^m -1)(q^{m} - q) \cdots (q^{m} - q^{\ell -1})}
{(q^{\ell} -1)(q^{\ell} - q) \cdots (q^{\ell} - q^{\ell -1})} \quad
{\rm and} \quad k = { {m} \choose {\ell} }.
\end{equation}
The minimum distance of $C(\ell ,m)$ is given by the following elegant formula 
due to Nogin \cite{N}: 
\begin{equation}
\label{MinDistClm}
d\left( C(\ell ,m)\right) = q^{\delta},  \quad {\rm where} \quad 
\delta := \ell (m - \ell).
\end{equation}
In fact, Nogin \cite{N} also determined some of the higher weights 
of $C(\ell ,m)$. More precisely, he showed that for
$1\le r \le \max \{\ell, m - \ell \} + 1$, 
\begin{equation}
\label{HgrWtClm}
d_r\left( C(\ell ,m)\right) = q^{\delta} +q^{\delta -1} + \cdots 
+ q^{\delta - r+1} .
\end{equation}
Alternative proofs of \eqref{HgrWtClm} were given in
\cite{GL}, and in the same paper a generalization to Schubert codes was 
proposed. The Schubert codes are indexed by the elements of the set
$$
I(\ell , m):= \{ \a = (\a_1, \dots , \a_{\ell} )\in \Z^{\ell} : 
1 \le \a_1 < \cdots < \a_{\ell} \le m\}.
$$
Given any $\a \in I(\ell , m)$, the corresponding {\em Schubert code} is denoted by
$C_{\a}(\ell ,m)$, and it is the code obtained from the 
projective system defined by the Schubert variety $\Omega_{\a}$ in
$G_{\ell, m}$ with a nondegenerate embedding induced by the Pl\"ucker embedding. 
Recall that $\Omega_{\a}$ can be defined by 
$$
\Omega_{\a} = \{ W \in G_{\ell, m} : \dim (W \cap A_{\a_i} ) \ge i \ {\rm for} \
i=1, \dots , \ell \},
$$
where $A_j$ denotes the span of the first $j$ vectors in a 
fixed basis of $V$, for $1\le j \le m$. 
It was observed in \cite{GL} that the length $n_{\a}$ and
the dimension $k_{\a}$ of $C_{\a}(\ell ,m)$ are abstractly given by 
\begin{equation}
\label{AbsNaKa}
n_{\a} = \vert \Omega_{\a} (\Fq) \vert \quad {\rm and} \quad
k_{\a} = \vert \{ \b \in I(\ell , m) : \b \le \a \} \vert,
\end{equation}
where 
for $\b =(\b_1, \dots , \b_{\ell} )\in I(\ell , m) $, by 
$\b \le \a$ we mean that $\b_i \le \a_i \ {\rm for} \ i=1, \dots , \ell$.
It was shown in \cite{GL} that the minimum distance of $C_{\a}(\ell ,m)$ 
satisfies the inequality
$$
d(C_{\a}(\ell ,m)) \le q^{\delta_{\a}},  \quad {\rm where} \quad
\delta_{\a} : = \sum_{i=1}^{\ell} (\a_i - i) = \a_1 + \cdots + \a_{\ell} - 
\frac{\ell (\ell + 1) }{2}.
$$
Further, it was conjectured by the first named author that, in fact, 
the equality holds, i.e.,
\begin{equation}
\label{MDC}
d(C_{\a}(\ell ,m)) = q^{\delta_{\a}}. 
\end{equation}
We shall refer to \eqref{MDC} 
as the {\em minimum distance conjecture} (for Schubert codes). 
Note that if $\a = (m - \ell +1, \dots , m -1, m)$, then $ \Omega_{\a} = G_{\ell, m}$ and 
so in this case  \eqref{MDC} is an immediate consequence of \eqref{MinDistClm}. 

The minimum distance conjecture has been proved in the affirmative by  
Hao Chen \cite{HC} when $\ell = 2$. In fact, he proves the following. 
If $\ell = 2$ and $\a = (m - h - 1, m)$ [we can assume that $\a$ is of this form
without any loss of generality], then 
$
d(C_{\a}(2,m)) = q^{\delta_{\a}} = q^{2m - h - 4},
$ 
and moreover, 
\begin{eqnarray}
\label{eq:HCn}
n_{\a} &=& \frac{ (q^m -1) (q^{m -1} - 1) } { (q^{2} -1) (q -1) } - 
\sum_{j=1}^h  \sum_{i=1}^{j}  q^{2m - j - 2 - i}, \ {\rm and}  \\
\label{eq:HCk}
k_{\a} &=& \frac{m (m - 1) }{2} - \frac{h (h + 1) }{2}.
\end{eqnarray}

An alternative proof of the minimum distance conjecture, as well as the weight
distribution of codewords in the case $\ell = 2$, was obtained 
independently by Guerra and Vincenti \cite{GV}; 
in the same paper, they prove also the following lower bound for 
$d(C_{\a}(\ell ,m))$ in the general case:
\begin{equation}
\label{GVLB}
d(C_{\a}(\ell ,m)) \ge \frac{q^{\a_1} (q^{\a_2} -q^{\a_1}) \cdots (q^{\a_{\ell}}
-q^{\a_{\ell - 1} })}{q^{1+2+\dots +\ell}} \ge q^{\delta_{\a} - \ell}.
\end{equation}
In an earlier paper, Vincenti \cite{V}, partly in collaboration with Guerra, 
verified the minimum distance conjecture for the unique nontrivial Schubert 
variety in the Klein quadric $G_{2,4}$, 
namely $\Omega_{(2,4)}$, 
and obtained a  lower bound which is weaker than \eqref{GVLB}, and also 
proved the following formula\footnote{In fact, in \cite{V} and \cite{GV}, the 
Grassmannian and its Schubert subvarieties are viewed as families of 
projective subspaces of a projective space rather than linear subspaces of 
a vector space. The two viewpoints are, of course, equivalent. To get 
\eqref{GVlength} from \cite[Prop. 15]{V}, one has to set $\ell =d+1$, 
$\a_i = a_{i-1}+1$ and $k_i = \ell_{i-1}+1$ for $1\le i \le \ell$. A
similar substitution has to be made to get \eqref{GVLB} from 
\cite[Thm. 1.1]{GV}.} 
for the length 
of $C_{\a}(\ell ,m)$.
\begin{equation}
\label{GVlength}
n_{\a} = \vert \Omega_{\a} (\Fq) \vert = \sum_{(k_1, \dots , k_{\ell - 1})} 
\prod_{i=0}^{\ell - 1}  { {\a_{i+1} - \a_{i}} \brack {  k_{i+1} - k_{i} }}_q
q^{ (\a_{i} - k_{i})(k_{i+1} - k_{i})},
\end{equation} 
where the sum is over all $(\ell -1)$-tuples $( k_1, \dots , k_{\ell - 1})$ of 
integers with $i\le k_i \le \a_i$ and $k_i \le k_{i+1}$ 
for $1\le i \le \ell - 1$, and where,
by convention, $\a_0 = 0 = k_0$ and $k_{\ell} = \ell$.

We can now describe the contents of this paper. In Section \ref{sec1} below,
we give two formulae for the length $n_{\a}$ 
of $C_{\a}(\ell ,m)$. Of these, the first is very simple and is related to a
classical result about the Grassmannians. The other formula is
somewhat similar to \eqref{GVlength} even though it was obtained independently.
The latter formula may be a little more effective in actual computations. 
Next, in Section \ref{sec2}, we give a
determinantal formula for the dimension $k_{\a}$ of $C_{\a}(\ell ,m)$ and
show that in certain cases this determinant can be evaluated. Moreover, we
also give an alternative formula for $k_{\a}$ using the formulae for $n_{\a}$
obtained in the previous section.
Finally, in Section \ref{sec3}, we show that the minimum distance and some of 
the higher weights for the codes corresponding to Schubert divisors, i.e.,  
Schubert varieties of codimension 
one in the corresponding Grassmannians, can be easily obtained using the
results of \cite{GL} and \cite{N}. This shows, in particular,
that the minimum distance conjecture  is true 
for all Schubert divisors such as, for instance, the unique nontrivial 
Schubert variety in the Klein quadric.

As a byproduct of the results in this paper, we see that
$n_{\a}$ can be expressed in three distinct ways and $k_{\a}$ in two.
This yields curious combinatorial identities, which may not be easy to prove directly.

Some of the main results of this paper, namely, Theorems 
\ref{Na2}, \ref{thm:ka} and \ref{DetEval}, were presented during a talk by the 
first author at the Conference on Arithmetic,
Geometry and Coding Theory (AGCT-8) held at CIRM, Luminy in May 2001. 
The article \cite{fpsac}, written for FPSAC-2003, gives an overview
(without proofs) of the results in this paper, and
it may be referred to for a more leisurely introduction to this paper.

We end this introduction with the following comment.
The Grassmannian is a special instance of homogeneous spaces
of the form $G/P$ where $G$ is a semisimple algebraic group and $P$ a parabolic
subgroup. Moreover, Schubert varieties also admit a generalization in this context.
Thus it was indicated in \cite{GL} that the Grassmann and Schubert codes can also
be introduced in a much more general setting. It turns out, in fact, that the 
construction of such general codes was already proposed in the binary case by
Wolper in an unpublished paper \cite{W}.
The general case, however, needs to be better understood and can be a source of
numerous interesting problems.

\section{Length of Schubert codes}
\label{sec1}

Fix integers $\ell, m$ with $1\le \ell \le m$. Let $I(\ell ,m)$ be
the indexing set with the partial order $\le $ defined in the previous 
section. For $\b =(\b_1, \dots , \b_{\ell} )\in I(\ell , m) $, let 
$$
\delta_{\b} : = \sum_{i=1}^{\ell} (\b_i - i) = \b_1 + \cdots + \b_{\ell} - 
\frac{\ell (\ell + 1) }{2}.
$$
Finally, fix some $\a \in I(\ell ,m)$ and let $C_{\a}(\ell ,m)$ be the 
corresponding Schubert code.

Quite possibly, the simplest formula for the length $n_{\a}$ of
$C_{\a}(\ell ,m)$ is the one given in the theorem below. This formula 
is an easy consequence of the well-known cellular decomposition of
the Grassmannian, which goes back to Ehresmann \cite{E}. However, it
doesn't seem easy to locate this formula in the literature, and thus, 
for the sake of completeness, we include here a sketch of the proof.

\begin{theorem}
\label{naCell}
The length $n_{\a}$ of
$C_{\a}(\ell ,m)$ or, in other words, the number of $\Fq$-rational points
of $\Omega_{\a}$, is given by 
\begin{equation}
\label{LengthCD}
n_{\a} = \sum_{\b \le \a} q^{ \delta_{\b} }
\end{equation}
where the sum is taken over all $\b \in I(\ell , m) $ satisfying $\b \le \a$. 
\end{theorem}

\begin{proof}
Consider, as in the previous section, the subspaces $A_j$ spanned by
the first $j$ basis vectors, for $1\le j \le m$. Given any $W \in
G_{\ell, m}$, the numbers $r_j = \dim W \cap A_j$ have the property\footnote{This 
follows, for example, because the kernel of the map $W \cap A_{j} \to \Fq$, 
mapping a vector to its $j$-th coordinate (with respect to the fixed basis 
of $V$), is $W \cap A_{j-1}$.} 
that $0\le r_j - r_{j-1} \le 1$ (where $r_0=0$, by convention), and, 
since $r_{m} = \ell$, there
are exactly $\ell$ indices where this difference is $1$. Thus there is a
unique $\b \in I(\ell , m) $ such that
$W$ is in
$$
C_{\beta} := \left\{ L \in  G_{\ell, m} : \dim (L \cap A_{\b_j} ) = j
\text{ and } \dim (L \cap A_{\b_j -1} ) = j -1 \text{ for }
1\le j \le \ell \right\}.
$$
Moreover,
for any $L \in C_{\beta}$, we have: $L \in \Omega_{\a} \Leftrightarrow
\b \le \a $. It follows that $\Omega_{\a}$ is the disjoint union of
$C_{\beta}$ as $\b$ varies over the elements of $I(\ell , m) $
satisfying $\b \le \a$. Now it suffices to observe that the subspaces
in $C_{\b}$ are in natural one-to-one correspondence with $\ell \times m$
matrices (over $\Fq$) with $1$ in the $(i,\b_i)$-th spot, and zeros to its
right as well as below, for $1 \le i \le \ell$.
\end{proof}

It may be argued that even though formula \eqref{LengthCD} is simple
and elegant, it may not be very effective in practice in view of the 
rather intricate summation involved. For example, if 
$\Omega_{\a}$ is the full Grassmannian $G_{\ell, m}$, then
\eqref{LengthCD} involves ${ {m} \choose {\ell} }$ summands, while the
closed form formula in \eqref{eq:1} given by the Gaussian binomial 
coefficient may be deemed preferable. For an arbitrary $\a\in I(\ell, m)$,
it is not easy to estimate the number of summands in \eqref{LengthCD}, as  
may be clear from the results of Section \ref{sec2}. 
With this in view, we shall now describe another formula for $n_{\a}$,
which is far from being elegant but may also be of some interest.
First, we need an elementary definition and a couple of preliminary lemmas.

By a {\em consecutive block} in an $\ell$-tuple
$\b =(\b_1, \dots , \b_{\ell} )\in I(\ell , m) $, we mean an ordered sequence
of the form $\b_i, \dots , \b_{j} $ where $1\le i \le j \le \ell$
and $\beta_{p+1} = \beta_p + 1$ for $i\le p < j$. For example, $3,4$
is a consecutive block in $(1,3,4,7)$ as well as in $(1,3,4,5)$ and 
in $(2,3,4,5)$. Note that any $\b \in I(\ell , m) $ always has $\ell$
consecutive blocks although it may often be regarded as having fewer
consecutive blocks.

\begin{lemma}
\label{ConOa}
Suppose $\a = (\a_1, \dots , \a_{\ell} )$ has $u+1$ consecutive blocks:
$$
\a = (\a_1, \dots , \a_{p_1}, \; \a_{p_1+1}, \dots , \a_{p_2}, \; \dots,
 \; \a_{p_{u-1}+1}, \dots , \a_{p_u}, \;
 \a_{p_{u}+1}, \dots , \a_{\ell} )
$$
so that $1\le p_1 < \cdots < p_u <  \ell$ and $\a_{p_i + 1}, \dots ,
\a_{p_{i+1}}$ are consecutive for $0 \le i \le u$, where by convention,
$p_0 =0$ and $p_{u+1} = \ell$. Then 
$$
\Omega_{\a} = \{ W \in G_{\ell, \a_{\ell} }: 
\dim (W \cap A_{\a_{p_i}}) \ge p_i \ {\rm for} \ i=1, \dots , u \}.
$$
\end{lemma}
 
\begin{proof}
As in the proof of Theorem \ref{naCell},  for any $W \in G_{\ell, \a_{\ell} }$, 
we have $\dim (W \cap A_{j-1}) \ge \dim (W \cap A_{j}) -1$ for $1\le j \le
m$.
Also, $\dim (W \cap A_{\a_{\ell}}) \ge \ell$ if and only if $W$ is a subspace
of $A_{\a_{\ell}}$. The desired result is now clear.
\end{proof}

Given any integers $a, b, s, t$, we define
$$
\lambda(a,b;s,t) =   \sum_{r=s}^{t} (-1)^{r - s} q^{{r - s}\choose {2}} 
{ {a - s} \brack {r - s}}_q { {b - r} \brack {t - r}}_q .
$$
Here, for any $u, v \in \Z$, the Gaussian binomial coefficient 
${ {u} \brack {v} }_q$ is defined
as in \eqref{eq:1} when $0\le v \le u$, and 
$0$ otherwise. 
Thus, 
if $a=s=0$, then $\lambda(a,b;s,t) = 
{ {b } \brack {t}}_q$. 

\begin{lemma}
\label{Mobius}
Let  $B$ be a $b$-dimensional vector space over $\Fq$ and $G_{t,b} =
G_t(B)$ denote the Grassmannian of $t$-dimensional subspaces of $B$.
Now suppose $A$ is any subspace of $B$ and $S$ is any subspace of
$A$, and we let $a = \dim A$ and $s=\dim S$. Then 
$$
\vert \{T\in G_t(B) : T \cap A = S \} \vert = \lambda(a,b;s,t) .
$$
\end{lemma}
 
\begin{proof}
Let $\LA$ be the poset of all subspaces of $A$ with the partial order
given by inclusion. Define functions $f, g: \LA \to \NN$ by 
$$
f(S) = \vert \{T\in G_t(B) : T \cap A = S \} \vert \quad
{\rm and}\quad g(S) = \vert \{T\in G_t(B) : T \cap A \supseteq S\}\vert .
$$
It is clear that for any $S\in \LA$ with $\dim S = s$, we have
$$
g(S) = {\mathop{ \sum_{R\in\LA} }_{R \supseteq S} }f(R) .
$$
On the other hand, for any $S$ as above, we clearly have 
\begin{equation}
\label{qting}
g(S) = \vert \{T\in G_t(B) : T \supseteq S \}\vert  = | G_{t-s}(B/S)|
= { {b - s} \brack {t - s} }_q.
\end{equation}
Hence, by M\"obius inversion applied to the poset $\LA$ and 
the well-known formula for the M\"obius function 
of $\LA$ (cf. \cite[Ch. 3]{StB}), we obtain
$$
f(S) = {\mathop{ \sum_{R\in\LA} }_{R \supseteq S} } \mu(S, R) g(R) 
= {\mathop{ \sum_{R\in\LA} }_{R \supseteq S} }
(-1)^{\dim R - \dim S} q^{ {{\dim R - \dim S}\choose{2}}} { {b - r}
\brack {t - r}}_q.
$$
Since the terms in the last 
summation depend only on the dimension of the varying subspace $R$,
we may write it as
$$
 \sum_{r=s}^{a} \vert \{R\in\LA : R \supseteq S \text{ and } \dim
R = r\}\vert  (-1)^{r - s} q^{{r - s}\choose {2}} { {b - r} \brack {t -
r}}_q .
$$
As in \eqref{qting}, the cardinality of the set appearing in the above summand 
is readily seen to be ${ {a - s} \brack {r - s}}_q $. This yields the
desired equality.
\end{proof}

\begin{theorem}
\label{Na2}
Let $u$ and $p_1, \dots , p_u$ be as in Lemma \ref{ConOa}. Then 
the length 
$n_{\a}$ 
of 
the Schubert code 
$C_{\a}(\ell ,m)$ is given by 
\begin{equation}
\label{eq:n}
n_{\a} = \sum_{s_1 = p_1}^{\a_{p_1}} \sum_{s_2 = p_2}^{\a_{p_2}} \cdots
\sum_{s_u = p_u}^{\a_{p_u}}
\; \prod_{i=0}^u \lambda(\a_{p_i} , \a_{p_{i+1}}; s_i, s_{i+1})
\end{equation}
where, by convention, $s_0 = p_0 = 0$ and $s_{u+1} =p_{u+1} = \ell$.
\end{theorem}

\begin{proof}
We use induction on $u$. If $u=0$, i.e., if $\a_1, \dots , \a_{\ell}$ are
consecutive, then $\Omega_{\a} = G_{\ell, \a_{\ell} }$, and so we know
that $n_{\a} = { { \a_{\ell}} \brack {\ell}}_q = \lambda(0,
\a_{\ell}; 0, \ell)$. Now suppose that $u\ge 1$ and the result holds for all
smaller values of $u$. Then, by Lemma \ref{ConOa}, we see that 
$$
\Omega_{\a} = \coprod_{S} \{T\in G_{\ell, \a_{\ell} } : T \cap A_{\a_{p_u}}
= S\}
$$
where the disjoint union is taken over the set, say $\Lambda_u$, of all
subspaces $S$ of
$A_{\a_{p_u}}$ satisfying $\dim S \ge u$ and $\dim S \cap A_{\a_{p_i}}
\ge p_i$ for $1 \le i \le u-1$. Hence, by Lemma \ref{Mobius}, 
$$
n_{\a} = \vert \Omega_{\a} (\Fq) \vert = 
\sum_{s=p_u}^{\a_{p_u}} 
\vert \{S \in \Lambda_u : \dim S = s \}\vert  \lambda(\a_{p_u} ,
\a_{\ell}; s, \ell).
$$
But for any $s$ with $p_u \le s\le \alpha_{p_u}$, the set of $s$-dimensional 
subspaces in $\Lambda_u$
is precisely the Schubert variety in $G_{s,\a_{p_u}}$ corresponding
to the tuple $(\a_1, \dots , \a_{p_u})$ with $u$ consecutive blocks. Hence
the induction hypothesis applies.
\end{proof} 

\begin{remark}
\label{chenNa}
In the case $\ell =2$, we obviously have $u \le 1$, and the formula given
above becomes somewhat simpler. It is not difficult to verify that this
agrees with the formula \eqref{eq:HCn} of Hao Chen \cite{HC}.
\end{remark} 

\begin{remark}
\label{combid1}
As a consequence of the results in this section, we obtain a purely combinatorial
identity which equates the right hand sides of \eqref{GVlength}, \eqref{LengthCD}
and \eqref{eq:n}. It would be an intriguing problem to prove this without invoking
Schubert varieties.
\end{remark} 

\section{Dimension of Schubert codes}
\label{sec2}

Let the notation be as in the beginning of the previous section. Our aim
is to give an explicit formula for the dimension $k_{\a}$  
of the Schubert code $C_{\a}(\ell ,m)$. As in the case of Theorem \ref{naCell},
it suffices to appeal to another classical fact about Schubert varieties in
Grassmannians, namely, the postulation formula due to Hodge \cite{Ho}. 
For our purpose, we use a slightly simpler description of Hodge's formula,
which (together with an alternative proof) is given in 
\cite{Gh}.

\begin{theorem}
\label{thm:ka}
The dimension 
$k_{\a}$ of the Schubert code $C_{\a}(\ell ,m)$ equals the
determinant of the $\ell \times \ell$ matrix whose $(i,j)$th entry
is ${ {\a_j - j +1} \choose {i - j + 1}}$, i.e., 
\begin{equation}
\label{eq:k}
k_{\a} = 
\left| \begin{array}{ccccc} 
{{\a_1} \choose {1}} & 1 &0 & \dots & 0 \\
{{\a_1} \choose {2}} &{{\a_2 - 1} \choose {1}} & 1 & \dots & 0 \\
\vdots & &  &   & \vdots \\{
{\a_1} \choose {\ell}} &{{\a_2 - 1} \choose {\ell - 1}} & 
{{\a_3 - 2} \choose {\ell - 2}} & \dots & {{\a_{\ell} - \ell + 1}\choose {1}} \\
\end{array}
\right|.
\end{equation}
\end{theorem}

\begin{proof}
Recall the abstract description in \eqref{AbsNaKa} 
for the dimension $k_{\a}$ of $C_{\a}(\ell ,m)$: 
$$
k_{\a} = \vert \{ \b \in I(\ell , m) : \b \le \a \} \vert.
$$
By Hodge Basis Theorem (cf. \cite[Thm. 1]{Gh}), we know that a 
vector space basis
for the $t$-th component, say $R_t$, of the homogeneous coordinate ring of 
$\Omega_{\a}$
is indexed by the $t$-tuples $\left(\b^{(1)}, \dots, \b^{(t)} \right)$ of
elements of $I(\ell, m)$ satisfying $\b^{(1)}\le \dots \le \b^{(t)} \le \a$.  
The postulation formula of Hodge gives the Hilbert function 
$h(t) = \dim R_t$ ($t\in \NN$) of this ring. Now, using \cite[Lemma 7]{Gh}, we
may write 
$$
h(t) = \det_{1\le i, j \le \ell }
 \left( { {t+ \a_j - j } \choose {t + i - j } } \right) \quad \text{ for } t\in \NN.
$$
By putting $t=1$, we get the desired result. 
\end{proof} 

\begin{remark}
\label{chenKa}
In the case $\ell =2$, we obviously have 
$$
k_{\a} = \a_1 (\a_2 - 1) - {{\a_1}\choose{2}} = \frac{\a_1 (2\a_2 - \a_1 -1)}{2}$$
and if we write $\a = (m-h-1, m)$, then we retrieve 
the formula \eqref{eq:HCk} of Chen \cite{HC}.
\end{remark} 

The determinant in \eqref{eq:k}
is not easy to evaluate in general. For example, none of the recipes in the
rather comprehensive compendium of Krattenthaler \cite{Kr} seem to be 
applicable. The following Proposition shows, however, that in a special case 
a simpler formula can be obtained. 

\begin{theorem}
\label{DetEval}
Suppose $\a_1, \dots , \a_{\ell}$ are
in an arithmetic progression, i.e., there are $c,d\in \Z$ such that
$\a_i = c(i-1)+d$ for $i=1,\dots , \ell$.
Let $\a_{\ell +1} = c\ell + d = \ell \a_2 + (1 - \ell ) \a_1$.
Then 
$$
k_{\a} = {\frac{\a_1}{\ell !} } \prod_{i=1}^{\ell - 1}(\a_{\ell +1} - i) 
= {\frac{\a_1}{\a_{\ell +1}}} {{\a_{\ell +1} } \choose { \ell }}.
$$
\end{theorem}

\begin{proof}
If 
$\a_i = c(i-1)+d$ for $i=1,\dots , \ell$, then the $(i,j)$-th entry of the
transpose of the $\ell \times \ell$ matrix  in \eqref{eq:k} can be written as 
$$
 { {c(i-1) + d - i +1} \choose {j - i + 1}} =
{ {BL_i +A} \choose {L_i  + j}} , \quad \text{where} \ 
B = 1-c, \ L_i = 1-i \text{ and } A=d.
$$
Now we use formula (3.13) in \cite[Thm. 26]{Kr}, which says that 
for an $\ell \times \ell$ matrix whose  $(i,j)$-th entry of the 
form 
${ {BL_i +A} \choose {L_i  + j}}$  [where $A,B$ can be indeterminates 
and the $L_i$'s are integers], the determinant is given by 
$$
{\frac{\prod_{1\le i<j\le \ell } (L_i - L_j)}{ \prod_{i=1}^\ell  (L_i+\ell )!} }
 \prod_{i=1}^\ell  {\frac{(BL_i+A)!}{\left( (B-1)L_i+A - 1\right)! } }
\prod_{i=1}^\ell  (A - Bi + 1)_{i-1},
$$
where in the last product we used the shifted factorial notation, viz.,
$(a)_0=1$ and  
$(a)_t = a(a+1) \cdots (a+t-1)$, for $t\ge 1$.
Substituting $B = 1-c, \ L_i = 1-i \text{ and } A=d$ and making elementary 
simplifications, we obtain the desired formula.
\end{proof} 

\begin{remark}
\label{conKa}
The simplest case, where the above Proposition is applicable is when 
$\a_1, \dots , \a_{\ell}$ are consecutive, i.e., $c=1$ and $\a_i = d+i-1$. 
Notice that in this
case, the formula for $k_{\a}$ reduces to ${{ d + \ell -1} \choose {\ell }}$.
Of course, this is not surprising since $\Omega_{\a}$ is nothing but the
smaller Grassmannian $G_{\ell, d + \ell -1}$ in this case. Thus, in this case 
we also have simpler formulae for $n_{\a}$ and $\delta_{\a}$ 
and the minimum distance conjecture is true. However, even in
this simplest case, the evaluation of the determinant in \eqref{eq:k} does
not seem obvious. Indeed, here it becomes an instance of the Ostrowski 
determinant
$\det \left( { {d}\choose{k_i -j} }\right)$ if we take $k_i = i+1$. 
A formula for such a determinant and the result that it is positive 
for increasing $\{k_i\}$ was obtained by Ostrowski \cite{Os} in 1964.
The case when $\{k_i\}$ are consecutive seems to go back to Zeipel in 1865 
(cf. \cite[Vol. 3, pp. 448-454]{Mu}).
\end{remark} 

An alternative formula for the dimension $k_{\a}$ of $C_{\a}(\ell, m)$ can be
derived using results of the previous section. To this end, we begin 
by observing that the dimension $k$ of the $q$-ary Grassmann code $C(\ell, m)$
doesn't depend on $q$, and bears the following relation to the length 
$n = n(q)$ of $C(\ell, m)$:
\begin{equation}
\label{limit:n}
\lim_{q\to 1} n(q) = k \quad \text{or, in other words, } \quad 
\lim_{q\to 1} { {m}\brack{\ell}}_q = { {m} \choose{\ell}}.
\end{equation}
Much has been written on this limiting formula in combinatorics literature. 
For example, a colourful, albeit mathematically incorrect, way
to state it would be to say that the (lattice of) subsets of an $m$-set is the
same as the (lattice of) subspaces of an $m$-dimensional vector space over the
field of one element!
In the proposition below, we observe that a similar relation holds in the case
of Schubert codes, and, then, use this relation to obtain the said alternative
formula for $k_{\a}$.

\begin{proposition}
\label{Limit}
The dimension $k_{\a}$ of the $q$-ary Schubert code $C_{\a}(\ell, m)$ is
independent of $q$ and is related to the length $n_{\a} = n_{\a}(q)$ of
$C_{\a}(\ell, m)$ by the formula
\begin{equation}
\label{limit:na}
\lim_{q\to 1} n_{\a}(q) = k_{\a}.
\end{equation}
Consequently, if $u$ and $p_1, \dots , p_u$ be are as in Lemma \ref{ConOa}, 
then 
\begin{equation}
\label{limit:ka1}
k_{\a} = \sum_{s_1 = p_1}^{\a_{p_1}} \sum_{s_2 = p_2}^{\a_{p_2}} \cdots
\sum_{s_u = p_u}^{\a_{p_u}}
\; \prod_{i=0}^u { {\a_{p_{i+1}} - \a_{p_i}}\choose{s_{i+1} - s_i}},
\end{equation}
where, by convention, $s_0 = p_0 = 0$ and $s_{u+1} =p_{u+1} = \ell$.
\end{proposition}

\begin{proof}
The limiting formula \eqref{limit:na} follows from the abstract description
in \eqref{AbsNaKa} of $k_{\a}$ and Theorem \ref{naCell}. Further, 
\eqref{limit:ka1} will follow from Theorem \ref{Na2} if we show that 
for any integer parameters $a,b,s,t$, we have 
$$
\lim_{q\to 1} \lambda(a,b;s,t) = { {b-a}\choose{t-s} }.
$$
But, in view of \eqref{limit:n}, this is equivalent to proving the 
binomial identity
$$
\sum_{j\ge 0} (-1)^j  { {a-s}\choose{j} } { {b-s-j}\choose{t-s-j} } 
= { {b-a}\choose{t-s} }.
$$
This identity is trivial if $t<s$, and if $t\ge s$, it follows easily if, 
after expanding by the binomial theorem, we compare the coefficients of $X^{t-s}$ in
the identity
$$
(1-X)^{a-s} (1-X)^{t-b-1} = (1-X)^{a-b+t-s-1} 
$$
and observe that for any integers $M$ and $N$, we have 
$
 { {-N-1}\choose{M} } = (-1)^M { {N+M}\choose{M} } .
$
\end{proof}

\begin{remark}
\label{combid2}
As a consequence of the results in this section, we obtain a purely combinatorial
identity which equates the right hand sides of \eqref{eq:k}
and \eqref{limit:ka1}. It would be an intriguing problem to prove this without invoking
Schubert codes.
\end{remark} 

While one would like to construct codes having both the {\em rate} $k/n$ and
the {\em relative distance} $d/n$ as close to $1$ as possible, 
the two requirements are in conflict with each other. 
For Schubert codes, this conflict manifests itself in a peculiar way:

\begin{corollary}
\label{Conflict}
Let $R = R(q)$ and $\Delta = \Delta (q)$ denote, respectively, the rate 
and the relative distance of the $q$-ary Schubert code  
$C_{\a}(\ell, m)$. Then, we have
$$
\lim_{q\to 1} R(q) = 1 \quad {\rm and } \quad 
\lim_{q\to \infty} \Delta (q)= 1.
$$ 
\end{corollary}

\begin{proof}
The limiting formula for the rate is immediate from Proposition \ref{Limit}.
As for the relative distance, it suffices to observe that using 
Theorem \ref{naCell}, we have 
$$
\lim_{q\to \infty} \frac{U_{\alpha}(q) }{n_{\a} (q)}= 1 \quad {\rm and } \quad 
\lim_{q\to \infty} \frac{L_{\alpha}(q) }{n_{\a} (q)}= 1,
$$
where $U_{\a}(q):= q^{\delta_{\a}}$ denotes the upper bound 
(cf. \cite[Prop. 4]{GL}) for the 
minimum distance of $ C_{\a}(\ell ,m)$, while
$L_{\a}(q)$ denotes the lower bound given by  \eqref{GVLB}.
\end{proof}

\section{Minimum Distance Conjecture for Schubert Divisors}
\label{sec3}


The notation in this section will be as in the Introduction and at the beginning
of Section \ref{sec1}. To avoid trivialities, we may tacitly assume 
that $1 < \ell < m$. Further, we let
$$
\theta:=(m-\ell+1, m-\ell+2, \dots ,m) \quad  \text{and} \quad
\eta:=(m-\ell, m-\ell+2, \dots ,m). 
$$
Note that with respect to the partial order $\le$, defined in the Introduction,
$\theta$ is the unique maximal element of $I(\ell, m)$ 
whereas $\eta$ the unique submaximal element.
Moreover, by \eqref{AbsNaKa}, we have 
$$
k_{\theta} = k := { {m} \choose {\ell} } \text{ and }
k_{\eta} = k -1 ; \text{ also }
\delta_{\theta} = \delta := \ell(m - \ell)  \text{ and } 
\delta_{\eta} = \delta - 1 . 
$$
Thus, in view of Theorem \ref{naCell}, we have 
\begin{equation}
\label{ntAne}
n_{\theta} = | \Omega_{\theta} | = |G_{\ell, m} | = 
 { {m} \brack {\ell}}_q \quad {\rm and} \quad 
n_{\eta} = | \Omega_{\eta} |  = 
 { {m} \brack {\ell}}_q  - q^{\delta} .
\end{equation}
Indeed, $\Omega_{\theta}$ is the full Grassmannian 
$G_{\ell, m}$, whereas $\Omega_{\eta}$
is the unique subvariety of $G_{\ell, m}$ of codimension one, which is often
referred to as the {\em Schubert divisor} in $G_{\ell, m}$. 

\begin{theorem}
\label{thm:dra}
If $\eta:=(m-\ell, m-\ell+2, \dots ,m)$ so that
$\delta_{\eta}=\delta -1 $, 
then  
\begin{equation}
\label{HgrWtOeta}
d_r\left( C_{\eta}(\ell ,m)\right) = q^{\delta -1 } +q^{\delta -2} + \cdots + q^{\delta - r} 
\quad {\rm for} \quad 1\le r \le \max \{\ell, m - \ell \}.
\end{equation}
In particular,
$d_1\left( C_{\eta}(\ell ,m)\right) = q^{\delta_{\eta} }$, and so the 
minimum distance conjecture is valid in this case.
\end{theorem}

\begin{proof}
Let $r$ be a positive integer and 
$H_{\theta} = \left\{ p=(p_{\beta})\in \PP^{k-1} = \PP(\wedge^{\ell}V) : 
p_{\theta} = 0 \right\}$ be the hyperplane given by the vanishing of the 
Pl\"ucker coordinate corresponding to $\theta$. Note that 
$\Omega_{\eta}=G_{\ell, m}\cap H_{\theta}$. Now, if $\Pi$ is a linear 
subspace of $\PP^{k_{\eta}-1} = \PP\left(H_{\theta} \right)$ of codimension $r$, 
then as a linear subspace of $\PP^{k-1}$, it is of codimension $r+1$. Therefore, 
$$
| \Omega_{\eta} \cap \Pi | = | G_{\ell, m}\cap H_{\theta} \cap {\Pi}|
 = | G_{\ell, m}\cap {\Pi}| \le 
  | G_{\ell, m} | - d_{r+1}\left( C(\ell, m) \right).
$$
Hence, in view of \eqref{ntAne}, if $r \le \max \{\ell, m -\ell \}$, 
then by \eqref{HgrWtClm}, we see that
$$
d_r\left( C_{\eta}(\ell ,m)\right) = 
| \Omega_{\eta} | - 
\max_{\codim \Pi = r} | \Omega_{\eta} \cap \Pi | \ge 
| \Omega_{\eta} | - | G_{\ell, m}| +  
q^{\delta } +q^{\delta -1} + \cdots + q^{\delta - r}.
$$
Thus, to complete the proof it suffices to exhibit a codimension $r$ 
linear subspace $\Pi$ of  
$\PP^{k_{\eta}-1} = \PP\left(H_{\theta} \right)$ such that
$| \Omega_{\eta} \cap \Pi | = | \Omega_{\eta} | - 
\left(q^{\delta -1 } +q^{\delta -2} + \cdots + q^{\delta - r}\right)$.
To this end, we use the notion 
of a close family introduced in \cite{GL} and \cite{FFA}, 
and some results from \cite{GL}.

First, suppose $m-\ell \ge \ell$ so that $r \le m - \ell$. 
Now let 
$$
\alpha^{(j)} = (m-\ell+2-j, m-\ell+2, m-\ell+3,\dots ,m), 
\quad {\rm for} \  j =1, \dots , r +1, 
$$
and let
$\Lambda = \left\{ \alpha^{(1)}, \dots , \alpha^{(r +1)} \right\}$. 
Then $\Lambda$ is a subset of $I(\ell, m)$ and a {\em close family}\footnote{Two 
elements $\alpha = (\a_1, \dots , \a_{\ell})$ and 
$\beta= (\beta_1, \dots , \beta_{\ell})$ in $I(\ell,m)$ are said to be {\em close}
if they differ in a single coordinate, that is, 
$|\left\{\a_1, \dots , \a_{\ell} \right\} \cap 
\left\{\beta_1, \dots , \beta_{\ell}\right\}| = \ell -1$. A subset of $I(\ell, m)$ 
is called a {\em close family} if any two distinct elements in it are close.}, 
in the sense of \cite[p. 126]{GL}. Note that 
$\alpha^{(1)} = \theta$ and $\alpha^{(2)} = \eta$.
Thus if $\Pi$ denotes the linear subspace of 
$\PP^{k_{\eta}-1} = \PP\left(H_{\theta} \right)$ defined by the
vanishing of the Pl\"ucker coordinates corresponding to 
$\alpha^{(2)}, \dots , \alpha^{(r + 1)}$, and 
$\Pi'$ denotes the linear subspace of 
$\PP^{k-1}$ defined by the
vanishing of the Pl\"ucker coordinates corresponding to 
$\alpha^{(1)}, \dots , \alpha^{(r + 1)}$, then 
$\codim {\Pi}' = r + 1 $, and 
using \cite[Prop. 1]{GL}, we obtain
$$
| \Omega_{\eta} \cap \Pi | = 
 | G_{\ell, m}\cap {\Pi}'| = 
 { {m} \brack {\ell}}_q  - q^{\delta} 
- q^{\delta -1 } - \cdots - q^{\delta - r}.
$$
Thus, in view of  \eqref{ntAne}, it follows that $\Pi$ is a subspace of 
$\PP^{k_{\eta}-1} = \PP\left(H_{\theta} \right)$ of codimension $r$ 
with the desired property.

On the other hand, suppose $\ell \ge m-\ell$. Then we let 
$$
\alpha^{(j)} = (m-\ell, m-\ell+1, \dots , \widehat{m-\ell+j-1}, \dots , m ),
\quad {\rm for} \  j =1, \dots , r +1, 
$$ 
where $\widehat{m-\ell+j-1}$ indicates that the element 
${m-\ell+j-1}$ is to be removed. Once again, for $r \le \ell$, 
$\Lambda = \left\{ \alpha^{(1)}, \dots , \alpha^{(r +1)} \right\}$
is a subset of $I(\ell, m)$ and a close family with 
$\alpha^{(1)} = \theta$.
Hence we can proceed as before and apply \cite[Prop. 1]{GL} to obtain
the desired formula for $d_r\left( C_{\eta}(\ell ,m)\right)$.
\end{proof} 

\begin{remark}
\label{inductived1}

An obvious analogue of the inductive argument in the above proof seems to
fail for Schubert subvarieties of codimension $2$ or more. 
For example, in $G_{3,6}$
the subvariety $\Omega_{\a}$ corresponding to $\a = (3,4,6)$ is of codimension
$2$. However, $\Omega_{\a}$  is not the intersection of $G_{3,6}$ with two
Pl\"uker coordinate hyperplanes but with four of them [viz., those corresponding to
$(j, 5,6)$ for $1\le j\le 4$]. Thus, to determine $d_1(C_{\alpha}(3,6))$, we 
should know $d_5(C(3, 6))$. But we know $d_r(C(3,6))$ only for 
$r\le\max\{3, 6-3\} + 1 = 4$. The argument will, however, work for Schubert 
varieties of codimension $2$ in $G_{2,m}$ because one of these two varieties
will be a lower order Grassmannian while the other is a section by just  3
hyperplanes, and assuming, as we may, that $m>4$,  we can apply formula
\eqref{HgrWtClm} and some results from \cite{GL}. We leave the details to the
reader. In any case, we know from the work of Hao Chen \cite{HC} and 
Guerra-Vincenti \cite{GV} that the minimum distance conjecture 
is true when $\ell =2$. 
\end{remark}

\section*{Acknowledgments}
\label{secAck}
\begin{small}
A part of this work was done when the first author was visiting the Institut
de Math\'ematiques de Luminy (IML), Marseille during May-June 1999. Partial 
support for this visit from the Commission on Development and Exchanges (CDE) 
of the International Mathematical Union (IMU), and the warm hospitality of 
IML is gratefully acknowledged. Thanks are also due to
Gilles Lachaud and Serge Vl\u{a}du\c{t} for useful discussions, to 
Christian Krattenthaler
for some helpful correspondence concerning determinant evaluations, and to 
Lucio Guerra and Rita Vincenti for a week long visit in June 2002 to the 
University of Perugia during which the first author learned of their 
recent results. 
We are also grateful to a vicious landlady, whom we do not name here, and more 
so to a kind-hearted one, Mme C. Cerisier, whose actions led us to share an 
apartment in Marseille in June 1996, where this work began.
\end{small}

\end{document}